\documentclass[12pt]{amsart}
\usepackage{amsmath}
\usepackage{amssymb}
\usepackage{amsthm}
\usepackage[english]{babel}
\usepackage{xcolor}
\usepackage{cite}
\usepackage{bigfoot}
\usepackage{hyperref}

\theoremstyle{plain}
\newtheorem{thm}{Theorem}[section]
\newtheorem{lem}[thm]{Lemma}
\newtheorem{prop}[thm]{Proposition}
\newtheorem{cor}[thm]{Corollary}
\theoremstyle{definition}
\newtheorem{remark}[thm]{Remark}

\newcommand{\eps}{\varepsilon}
\newcommand{\Z}{\mathbb{Z}}
\newcommand{\R}{\mathbb{R}}
\newcommand{\AN}{\mathcal{A}_N}

\newcommand{\calE}{\mathcal{E}}
\newcommand{\calF}{\mathcal{F}}

\newcommand{\calH}{\mathcal{H}}
\newcommand{\calI}{\mathcal{I}}
\newcommand{\calK}{\mathcal{K}}
\newcommand{\calP}{\mathcal{P}}
\newcommand{\calR}{\mathcal{R}}

\newcommand{\calT}{\mathcal{T}}
\newcommand{\calW}{\mathcal{W}}
\newcommand{\p}{\partial}
\newcommand{\Ds}{(-\Delta)^{s}}
\newcommand{\Dn}{(-\gamma^2\Delta)^{-1}}
\newcommand{\Dnh}{(-\gamma^2\Delta)^{-\frac12}}
\newcommand{\textas}{\text{ as }}
\newcommand{\textin}{\text{ in }}
\newcommand{\texton}{\text{ on }}
\newcommand{\textif}{\text{ if }}
\newcommand{\textfor}{\text{ for }}

\DeclareMathOperator{\BV}{BV}

\newcommand{\set}[1]{\left\{#1\right\}}
\newcommand{\abs}[1]{\left|#1\right|}
\newcommand{\norm}[2][]{\left\|#2\right\|_{#1}}
\newcommand{\angles}[1]{\left\langle#1\right\rangle}

\begin{document}

\title[Lamellar phase solutions]{Lamellar phase solutions for diblock copolymers with nonlocal diffusions}

\author[H. Chan, M. Jamshid Nejad and J. Wei]{Hardy Chan, Masomeh Jamshid Nejad and Juncheng Wei}
\address[H. Chan]{Department of Mathematics, The University of British Columbia}
\email[H. Chan]{hardy@math.ubc.ca}
\address[M. Jamshid Nejad]{Department of Mathematics, The University of British Columbia}
\email[M. Jamshid Nejad]{masomeh@math.ubc.ca}
\address[J. Wei]{Department of Mathematics, The University of British Columbia}
\email[J. Wei]{jcwei@math.ubc.ca}


\begin{abstract}
For a diblock copolymer with total chain length $\gamma>0$ and mass ratio $m\in(-1,1)$, we consider the problem of minimizing the doubly nonlocal free energy
\[\mathcal{E}_{\eps}(u)
=\mathcal{H}(u)
+\frac{1}{\varepsilon^{2s}}
    \int_{\Omega}W(u)\,dx
+\frac{1}{2}\int_{\Omega}
    \abs{(-\gamma^{2}\Delta)^{-\frac{1}{2}}(u-m)}^2\,dx
\]
in a domain $\Omega$, where $\mathcal{H}(u)$ is a fractional $H^s$-norm with $s\in(0,\frac12)$, and $W$ is a double-well potential. This arises in the study of micro-phase separation phenomena for diblock copolymers with nonlocal diffusions.

On the unit interval, we identify the $\Gamma$-limit as $\eps\to0^+$, and also find explicit isolated local minimizers associated the lamellar morphology phase in the case $m=0$, provided that the chain is sufficiently short or the nonlocal interaction is sufficiently strong (i.e. as $s\to0^+$). We stress that such extra condition is new for the nonlocal case and is not present in the classical model. The proof, while elementary, requires a careful analysis of the nonlocal integrals.
\end{abstract}


\maketitle

\section{Introduction}

\subsection{Diblock copolymers}

In this article, we study the microphase separation phenomenon of diblock copolymers under the effect of nonlocal diffusions. Originally, the model, with ordinary diffusion, was introduced in Bahiana and Oono \cite{Bahiana-Oono90} and Ohta and Kawasaki \cite{Ohta-Kawasaki86}.

For a domain $\Omega\subset\R$ and $u\in L^2(\Omega)$, consider the following free energy
\begin{equation}\label{eq:Eeps}
\calE_\eps(u)=
\begin{cases}
\calH(u)+\dfrac{1}{\eps^{2s}}\calW(u)+\calK(u)
    &\textif{u}\in{X}\cap{H^s}(\Omega),\\
+\infty
    &\textif{u}\in{X}\setminus{H^s}(\Omega),
\end{cases}
\end{equation}
where
\begin{equation}\label{eq:calH}
\calH(u)=\dfrac14\int_{\Omega}\int_{\Omega}
    (u(x)-u(y))^2K_\Omega(x-y)\,dxdy,
\end{equation}
\begin{equation}\label{eq:calW}
\calW(u)=\int_{\Omega}W(u(x))\,dx
\end{equation}
\begin{equation}\label{eq:calK}
\calK(u)=
\dfrac12\int_{\Omega}\abs{\Dnh(u-m)(x)}^2\,dx
\end{equation}
and $K_\Omega$ is a fractional kernel of order $2s$ such that $\calH(u)$ is comparable to the $H^s$-norm. Here $m=\int_0^1u$ is the mass ratio. The admissible class $X$ of functions that takes values between $[-1,1]$ with a fixed mass ratio $m$ is introduced precisely in \eqref{eq:X} later in the Introduction. 

In a diblock copolymer, a linear-chain molecule consists of two sub-chains covalently grafted to each other. The sub-chains consist two monomer units which can be represented by $u=-1$ and $u=+1$, the global minima of the double-well potential $W$. On the one hand, the potential energy $\calW$ takes into account the repulsion of the two monomers. On the other hand, there are the long-range chemical bonding of the two monomer units. This part of the free energy is inversely proportional to the square root of the total chain length $\gamma$, as seen in the $\calK$ term. As a result, such segregation only leads to phase separation in the microscopic scale. Here $\eps$ is proportional to the thickness of interfaces. The interfacial energy is represented by $\calH$; it prevents the unnecessary forming of interfaces.

As we minimize the free energy $\calE_\eps$, let us first observe that the double-well potential $W$ prefers segregated monomers to a mixture. With a fixed mass ratio $m$ of the two monomers, we see two competing tendencies --- $\calH$ likes large blocks of monomers\footnote{
Indeed, from the formula, if $x$ and $y$ are points which are close to each other and lie on two sides of an interface, then $u(x)-u(y)$ is $O(1)$ while the singular fractional kernel at $x-y$ is very large.
} while $\calK$ favors rapid oscillations. The process of reaching a stable configuration is termed \emph{micro-separation}, and the patterns formed in micro-domains are known as \emph{morphology phases}.

The classical model has instead the energies
\[\dfrac{\eps}{2}\int_{\Omega}|\nabla{u}|^2\,dx
    +\dfrac{1}{\eps}\calW(u)\]
in place of $\calH(u)+\eps^{-2s}\calW(u)$ in the free energy \eqref{eq:Eeps}. This modification is made as we take into consideration the longer-range interactions around the interfaces, modeled by a nonlocal diffusion term. For the derivation of a similar model with nonlocal diffusion, the readers are referred to the appendix. Later in the introduction we will discuss the roles of the thickness parameter $\eps$ in the article in the context of $\Gamma$-convergence.

There is a great amount of literature concerning the classical model. With different ranges of the parameters $\eps$ and $\gamma$, various morphology phases of diblock copolymers are possible and they have been studied mathematically: the \emph{lamellar phase} by Ren--Wei \cite{Ren-Wei00, Ren-Wei03min, Ren-Wei03spectra, Ren-Wei05, Ren-Wei06}, Fife--Hilhorst \cite{Fife-Hilhorst01}, Choksi \cite{Choksi01}, Choksi--Ren \cite{Choksi-Ren05}, Chen--Oshita \cite{Chen-Oshita05}, Choksi--Sternberg \cite{Choksi-Sternberg07}, the \emph{cylindrical phase} by Ren--Wei \cite{Ren-Wei07single, Ren-Wei07many}, the \emph{spherical phase} by Glasner--Choksi \cite{Glasner-Choksi09}, Choksi--Peletier \cite{Choksi-Peletier10, Choksi-Peletier11} and the \emph{gyroid} and \emph{orthorhomic phases} by Teramoto--Nishiura \cite{Teramoto-Nishiura02, Teramoto-Nishiura10}. A detailed analysis of the phase diagram is given by Choksi--Peletier--Williams \cite{Choksi-Peletier-Williams09}. For more details of this model and the associated parabolic problem, the readers may consult \cite{Bahiana-Oono90, Mueller93, Nishiura-Ohnishi95, Escher-Nishiura02, Henry-Hilhorst-Nishiura03, Choksi-Ren03, Ren-Wei08}. As with the nonlocal diffusion, Dipierro--Novaga--Valdinoci \cite{Dipierro-Novaga-Valdinoci17} considered a nonlocal energy involving the fractional perimeter functional and established a rigidity result for critical points (not just minimizers) provided the volume is small in a certain sense.

\subsection{The two nonlocal terms}

In the recent decades there has been an explosive amount of literature concerning the fractional Sobolev space $W^{s,p}$ and the associated fractional Laplacian operator $\Ds$. They are useful in describing long-range interactions in physical systems including, for example, the Ising model in statistical mechanics, the Peierls--Nabarro model in dislocations in crystals, and the Benjamin--Ono equation in hydrodynamics. Recent progress and references can be found in \cite{Bucur-Valdinoci16}.

From a physical point of view, $-\Delta{u}$ describes the `usual' diffusion of certain particles with density $u$ due to random movements as it measures the deviation from the mean value in an infinitesimal neighborhood
. In contrast, fractional order operators like $\Ds{u}$, or more generally
\[Lu(x)=\int_{\Omega}(u(x)-u(y))K_\Omega(x-y)\,dy,\]
calculates the $K_\Omega$-weighted average in the whole domain. In such setting, particles do not only interact with its immediate neighbours, but are also allowed to influence others which are far apart.

To avoid complications, although the model can be easily generalized to a domain in $\R^d$, we impose periodic boundary conditions and work in one dimension only. When $\Omega=(0,1)$,
\begin{equation}\label{eq:K}
K_\Omega(x-y)=K(x-y)=C_{1,s}\sum_{n\in\Z}\dfrac{1}{|x-y-n|^{1+2s}}.
\end{equation}
For its derivation
, we refer to Roncal and Stinga \cite{Roncal-Stinga16}; see also \cite{DelaTorre-delPino-Gonzalez-Wei17, Ambrosio18}. Here $C_{1,s}=\frac{2^{2s}\Gamma(\frac{1+2s}{2})}{|\Gamma(-s)|\pi^{1/2}}$ is the normalization constant which also appears in the singular integral definition for the fractional Laplacian
\[\Ds{u}(x)=C_{1,s}\int_{\R}
    \dfrac{u(x)-u(y)}{|x-y|^{1+2s}}\,dy.\]
Indeed, for a Schwartz function, $\widehat{\Ds{u}}(\xi)=|\xi|^{2s}\hat{u}(\xi)$. (For a proof, see, for example, \cite{DiNezza-Palatucci-Valdinoci12}.)

\bigskip

Since the Laplacian
\begin{multline*}
-\gamma^{2}\Delta: \left\{{v\in H^{2}([0,1]): \int_{\Omega} v=0
\ \mbox{and} \ v \mbox{ is 1-perioidic} }\right\}
\\
\to \left\{{u\in L^{2}([0,1]):\int_{\Omega} u=0 \text{ and } u \text{ is 1-periodic} }\right\}
\end{multline*}
is an isomorphism, we may write $(-\gamma^{2}\Delta)^{-\frac{1}{2}}$ as the square root of the inverse of the Laplacian under the periodic boundary and zero average conditions, i.e. $v=\Dn{f}$ if and only if
\[\begin{cases}
-\gamma^2v''=f\quad\texton(0,1),\\
v(0)=v(1),\\
v'(0)=v'(1),\\
\displaystyle\int_0^1v=0.
\end{cases}\]
It is important to note that the presence of the nonlocal term $\calK$ 
actually gives us the local minima with 
a large number of
transitional layers which is in sharp contrast to local problems. In its absence, $\calE_\eps$ is reduced to the fractional Allen--Cahn energy whose unique (up to translation and reflection) global minimizer has a single layer \cite{Palatucci-Savin-Valdinoci13, Cabre-Sire15}; see also \cite{Gui-Zhao15}. Such energy is closely related to fractional minimal surfaces, first introduced by Caffarelli, Roquejoffre and Savin \cite{Caffarelli-Roquejoffre-Savin10}. Since then, the regularity, rigidity and qualitative behaviors of such surfaces have been widely studied. The interested readers are referred to the survey \cite{Valdinoci13}.

If $W\in C^{1}(\R)$, then a critical point $u$ of $\calE_\eps$ together with a $v$ and a $\lambda,$ solves the Euler--Lagrange equation
\[\begin{cases}
\Ds{u}+\dfrac{1}{\eps^{2s}}W'(u)+v=\lambda
    &\texton(0,1),\\
-\gamma^{2}v''=u-m
    &\texton(0,1),\\
u(0)=u(1),\quad
    v(0)=v(1),\\
u'(0)=u'(1),\quad
    v'(0)=v'(1),\\
\displaystyle\int_0^1u=0,\quad
    \displaystyle\int_0^1v=0.\\
\end{cases}\]

\subsection{$\Gamma$-convergence}

In this paper we show that as $\eps$ tends to $0$, $\calE_{\eps}$ converges to $\calE$, defined by
\begin{equation}\label{eq:E}
\calE(u)=
\begin{cases}
\calH(u)+\calK(u)
    &\textif{u}\in\BV([0,1],\set{-1,1}),\\
+\infty
    &\textif{u}\in{L^2}([0,1])\setminus\BV([0,1],\set{-1,1}).
\end{cases}
\end{equation}
Here $\BV([0,1],\set{-1,1})$ is the space of functions with bounded variations taking only the values $-1$ and $+1$. The convergence falls in the general theory of $\Gamma$-limit, on which there has been a considerable amount of literature. In \cite{DeGiorgi-Franzoni75}, De Giorgi and Franzoni introduced the $\Gamma$-convergence as `a notion of convergence for functionals, which tends to be as compatible as possible with the minimizing features of the energy, and whose limit is capable to capture essential features of the problem.' A notably important and relevant example was given by Modica and Mortola \cite{Modica-Mortola77}, who showed that the sequence of rescaled Allen--Cahn functionals
\[\calF_\eps(u)=
\begin{cases}
\displaystyle\int_\Omega\left(\dfrac{\eps}{2}|\nabla{u}|^2+\dfrac{1}{\eps}W(u)\right)\,dx
    &\textif{u}\in{H^1}(\Omega),\\
+\infty
    &\textif{u}\in{L^1}(\Omega)\setminus{H^1}(\Omega),
\end{cases}\]
$\Gamma$-converges as $\eps\to0^+$ to
\[\calF(u)=
\begin{cases}
c(W)\|Du\|(\Omega)
    &\textif{u}\in\BV(\Omega;\set{-1,1})\\
+\infty
    &\textif{u}\in{L^1}(\Omega)\setminus\BV(\Omega;\set{-1,1}),
\end{cases}\]
where $\Omega\subset\R^d$, $\|Du\|$ is the absolute value of the distributional derivative $Du$ as a finite measure, and
\[c(W)=\dfrac{1}{2}\displaystyle\int_{-1}^{1}\sqrt{2W(s)}\,ds.\]
Some results regarding nonlocal energies are \cite{Ambrosio-DePhilippis-Martinazzi11, Gonzalez09, Savin-Valdinoci12}. In particular, Savin and Valdinoci \cite{Savin-Valdinoci12} proved that the fractional energy $\calH(u)+\eps^{-2s}\calW(u)$ $\Gamma$-converges to the fractional perimeter functional\footnote{
For $s\in(0,\frac12)$, the fractional perimeter of a measurable set $E\subset\R^n$ in an open set $\Omega\subset\R^n$ is defined as the functional
\[
P_s(E,\Omega):=
\int_{E\cap\Omega}\int_{\R^n\setminus{E}}\dfrac{dxdy}{|x-y|^{n+2s}}
+\int_{E\setminus\Omega}\int_{\Omega\setminus{E}}\dfrac{dxdy}{|x-y|^{n+2s}},
\]
whenever the right hand side is finite.
} if $s\in(0,\frac12)$, while for $s\in[\frac12,1)$, a multiple of it $\Gamma$-converges to the classical perimeter. More precisely, we have
\[\begin{cases}
\calH(u)+\eps^{-2s}\calW(u)\xrightarrow{\Gamma}\calH(u)
	&\textfor s\in(0,\frac12),\\
\frac{1}{|\log\eps|}\calH(u)+\frac{1}{\eps|\log\eps|}\calW(u)\xrightarrow{\Gamma}c_{\frac12}(W)\norm{Du}(\Omega)
	&\textfor s=\frac12,\\
\eps^{2s-1}\calH(u)+\eps^{-1}\calW(u)\xrightarrow{\Gamma}c_s(W)\norm{Du}(\Omega)
	&\textfor s\in(\frac12,1),
\end{cases}\]
if $u\in\BV(\Omega;\set{-1,1})$, where $c_s$ is an constant depending on a one-dimensional profile. If $u\in\BV(\Omega;\set{-1,1})$, then the limit if $+\infty$.
Other classical examples of $\Gamma$-convergence are contained in \cite{Modica87, Braides02, DalMaso93, Tartar09} and the references therein.

With the $\calK$ term, Ren and Wei \cite{Ren-Wei00} found the $\Gamma$-limit of $\calF_\eps(u)+\calK(u)$ as $\calF(u)+\calK(u)$ in the ambient function space $L^2([0,1])$ with a fixed mass ratio, namely
\[\int_0^1u=m.\]
Note that their choice of $L^2$ as opposed to the classical $L^1$ is more natural in the presence of the $H^{-1}$ energy $\calK$. The local minima of the limiting problem are proved to be steps functions with evenly spaced jumps across $-1$ and $+1$, hence, in their neighborhood, the existence local minima of $\calF_\eps(u)+\calK(u)$ are also established for small $\eps$.

\subsection{Main results}

The aim of this paper is two-fold. First we identify the $\Gamma$-limit $\calE$ of the functional $\calE_{\eps}$. Then we find explicit local minimizers of $\calE$ (in fact also of $\calE_\eps$) that correspond to the lamellar morphology phase of a diblock copolymer.

While the general idea goes in parallel with \cite{Savin-Valdinoci12} and \cite{Ren-Wei00}, some comments are in order. Because of 
\cite[Theorems 1.2--1.3]{Savin-Valdinoci12}
, the regime $s\in[\frac{1}{2},1)$ is similar to the classical case $s=1$ and hence we will not study it here. Owing to the nonlocal terms, the problem makes sense only in the periodic boundary condition with the functions defined in the whole real line, as opposed to the Neumann boundary condition employed in \cite{Ren-Wei00}.

The difficulty lies exactly in the fractional norm $\calH$. In contrast to the local problem, the relevant system of equations is nonlocal and cannot be solved in the usual way. The novelty is therefore to find the candidate of the minimizer and to prove it with appropriate integral computations.


In order to state our first result, let us consider the function space
\begin{equation}\label{eq:X}
X=\set{u\in{L}^{\infty}([0,1]):u\text{ is 1-periodic}, \norm[L^{\infty}({[0,1]})]{u}\leq1,\text{ and } \int_{0}^{1}u\,dx=m},
\end{equation}
endowed with the topology of $L^{2}([0,1])$, i.e. we say that
\[u_\eps\xrightarrow{X}u\]
if
\[\int_{0}^{1}|u_\eps-u|^2\,dx\to0\quad\textas\eps\to0^+.\]

\begin{thm}\label{thm:Gamma}
Let $s\in(0,\frac{1}{2})$ and $\Omega=(0,1)$. For $\calE_\eps$ and $\calE$ defined in \eqref{eq:Eeps}--\eqref{eq:K} and \eqref{eq:E} respectively, $\calE_\eps$ $\Gamma$-converges to $\calE$ as $\eps\to0$, i.e. for any $u\in{X}$,
\begin{enumerate}
\item for any sequence $u_\eps\xrightarrow{X}u$, the liminf inequality holds, i.e.
    \[\calE(u)\leq\liminf_{\eps\to0^+}\calE_\eps(u_\eps);\]
\item there exists a recovery sequence $u_\eps\xrightarrow{X}u$ such that the limsup inequality holds, i.e.
    \[\limsup_{\eps\to0^+}\calE_\eps(u_\eps)\leq\calE(u).\]
\end{enumerate}
\end{thm}

Our second result concerns the explicit local minimizers of $\calE$. From now on we focus on the case $m=0$ where certain computations can be done explicitly. (We do not need this in Section \ref{sec:K}, though.) For an even integer $N\geq2$, let $\AN$ be the set of step functions of the form
\[u(x)=\sum_{k=0}^{N}(-1)^{k}\chi_{[x_k,x_{k+1}]}(x),\quad{x}\in(0,1),\]
with $0=x_0<x_1<x_2<\cdots<x_N<x_{N+1}=1$, such that
\begin{equation}\label{eq:intu}
\int_0^1u=0,\quad\text{ i.e. }
    2x_1-2x_2+\cdots-2x_N+1=0.
\end{equation}
An alternative and more useful expression in terms of the Heaviside step function $H$ is
\begin{equation}\label{eq:uH}
u(x)=1+2\sum_{k=1}^{N}(-1)^{k}H(x-x_k).
\end{equation}

Throughout the paper we write $U_N$ as the step function in $\AN$ such that the configuration points $(x_1,\dots,x_N)$ are \emph{equi-distributed}, meaning that
\[x_k=\dfrac{2k-1}{2N},\quad
    k=1,\dots,N.
\]

We have the following
\begin{thm}\label{thm:localmin}
Suppose $m=0$. For any even integer $N\geq2$, there exists an explicit $\gamma_0(N,s)>0$ such that for any $0<\gamma<\gamma_0(N,s)$, $U_N$ is a local minimizer of $\calE$ in $\AN$, in the sense that
\[D_{(x_1,\dots,x_N)}^{2}\calE(U_N)\geq0,\]
as a positive semi-definite matrix in the orthogonal complement \[E_{\frac{N}{2}}^{\perp}=\langle(1,-1,1,-1,\dots,1,-1)\rangle^{\perp}.\]
Moreover, the eigenvalues of the Hessian are given by
\[\begin{split}
\lambda_0&=0,\\
\lambda_{\ell}&=\dfrac{1}{\gamma^2{N}}
    \tan^2\left(\dfrac{\pi\ell}{N}\right)
        (1+O(\gamma^2)),
    \quad\textas\gamma\to0^+,
    \quad\textfor\ell\in\set{1,\dots,N}\setminus\set{\frac{N}{2}},\\
\lambda_{\frac{N}{2}}&=-\dfrac{4}{3\gamma^2N}(1+O(\gamma^2))
    \quad\textas\gamma\to0^+.
\end{split}\]
\end{thm}

There is no analogy of the condition $\gamma<\gamma_0(N,s)$ in the classical case. The standard perimeter functional does not change if one moves the interfaces of $U_N$ slightly. In our problem, however, the $s$-perimeter, the first term in $\calE$, changes if the interfaces of $U_N$ move. This induces a coarsening effect that penalizes multiple interfaces. This effect is overcome only if the second nonlocal term is sufficiently strong, i.e. $\gamma$ is sufficiently small.

\begin{remark}
From equation \eqref{eq:gamma0} in the proof, one may in fact take
\[\gamma_0(N,s)=\dfrac{1}{100\sqrt{s}N^{1+s}}\tan\left(\dfrac{\pi}{N}\right).\]
\end{remark}

\begin{remark}
We emphasize that $U_N$ is indeed a local minimizer, not just a saddle point. This is because any non-trivial variation in the direction $(1,-1,\dots)$ would violate the constraint \eqref{eq:intu}.
\end{remark}

\begin{remark}
Since $\calW\geq0$ and $\calW(u)=0$ for any $u\in\AN$, we see that $U_N$ are also local minimizers of $\calE_\eps$.
\end{remark}

\begin{remark}
In the more general case $m\neq0$, $x_k$ depends on $m$ as in \cite{Ren-Wei00} and one cannot expect a clean formula for the eigenvalues in the main order term like above. Nonetheless, we still expect the same result to be true, at least for $m$ close to $0$, via perturbative methods.

\end{remark}

In fact, we also see that $U_N$ are isolated local minimizers. This simply follows from \cite[Proposition 2.3]{Ren-Wei00} (now with the compact Sobolev embedding $H^s(\Omega)\hookrightarrow L^2(\Omega)$).

Finally, we compute the energy of the local minimizer $U_N$.

\begin{thm}\label{thm:EUN}
Let $m=0$. There exists a constant $C>0$ such that for any $N\geq1$,
\[C^{-1}\left(N^{2s}+\dfrac{1}{\gamma^2N^2}\right)
    \leq\calE(U_N)
        \leq{C}\left(N^{2s}+\dfrac{1}{\gamma^2N^2}\right).\]
\end{thm}


The paper is organized as follows. In Section \ref{sec:Gamma} we show the $\Gamma$-convergence of the free energy and prove Theorem \ref{thm:Gamma}. In Sections \ref{sec:H} and \ref{sec:K} we compute the derivatives of each term. Then in Section \ref{sec:localmin} we prove Theorem \ref{thm:localmin} by finding the explicit local minimizer. Next, in Section \ref{sec:v}, we obtain and explicit formula for $v$ that is useful for computing $\calK$ explicitly. Finally, an energy estimate concerning the growth of the fractional norm is contained in Section \ref{sec:growth}.

In the appendix we derive a diblock copolymer model that involves a fractional gradient. Note, however, that for the simplicity of the mathematical treatment, we decided to use the $H^{-1}$-norm which behaves similarly to the derived $H^{-s}$ norm and is anyway nonlocal.


\section{$\Gamma$-convergence and Existence of global minimizers}
\label{sec:Gamma}

Intuitively speaking, since for $s\in(0,\frac12)$, functions with jumps are allowed in $H^s$, we expect the $\Gamma$-limit
\begin{equation}
\calE(u)=\calH(u)+\calK(u)
\end{equation}
for $u=\chi_E-\chi_{E^c}\in\AN$, and $\calE(u)=+\infty$ otherwise. Since such convergence result is crucial for our purpose and the proof is short, for the sake of completeness, in this section we establish the $\Gamma$-convergence rigorously, following closely the argument in \cite{Savin-Valdinoci12}.

\begin{proof}[Proof of Theorem \ref{thm:Gamma}]
First we observe that if $u=\chi_{E}-\chi_{ E^c}$, then
\begin{equation}\label{eq:u=chi}
\calE_\eps(u)=\calE(u)=\calH(u)+\calK(u).
\end{equation}
To prove part (1), let $u_{\varepsilon}\xrightarrow{X}{u}$. This would be obvious if
\[\liminf_{\eps\to0^+}\calE_\eps(u_\eps)=+\infty.\]
Therefore, we suppose that
\[\liminf_{\eps\to0^+}\calE_\eps(u_\eps)=\ell<+\infty,\]
and that for a subsequence ${\eps_k}$ the above limit is actually attained. By passing to a further subsequence $u_{\eps_{k_j}}$, we may also assume that $u_{\eps_{k_j}}\to u$ almost everywhere. Now
\[\ell=\lim_{k\to+\infty}\calE_{\eps_k}(u_{\eps_k})
    =\lim_{j\to+\infty}\calE_{\eps_{k_j}}(u_{\eps_{k_j}})
    \geq\lim_{j\to+\infty}\dfrac{1}{\eps_{k_j}^{2s}}\int_{0}^{1}W(u_{\eps_{k_j}}(x))\,dx,\]
which implies that
\[\int_{0}^{1}W(u(x))\,dx
    =\lim_{j\to+\infty}\int_{0}^{1}W(u_{\eps_{k_j}}(x))\,dx
    =0.\]
This forces $u(x)\in\set{-1,+1}$ for almost every $x\in(0,1)$, meaning that $u=\chi_E-\chi_{E^c}\in\AN$ for a suitable set $E$. Since the energies $\calH$ and $\calK$ are lower semicontinuous, from \eqref{eq:u=chi} we have
\[\begin{split}
\calE(u)
&=\calH(u)+\calK(u)
    \\
&\leq\liminf_{\eps\to0^+}
    \left(\calH(u_\eps)+\calK(u_\eps)
    \right)\\
&=\liminf_{\eps\to0^+}\calE_\eps(u_\eps).
\end{split}\]
This proves (1).

Now we turn to the proof of (2). We may assume that $u=\chi_E-\chi_{E^c}\in\AN$, for otherwise the statement is vacuously true. In such case, we simply take $u_\eps=u$ and use \eqref{eq:u=chi} to conclude that
\[\calE(u)=\calE_\eps(u_\eps)
    \geq\limsup_{\eps\to0^+}\calE_\eps(u_\eps),\]
as desired.
\end{proof}

We recall the following adaptation \cite{Palatucci-Savin-Valdinoci13} of the classical Riesz--Fr\'{e}chet--Kolmogorov theorem on the compactness of a subset of $L^2([0,1])$ uniformly bounded in an $H^s$ norm.

\begin{lem}
Let $\calT\subset{L}^2([0,1])$ be bounded such that
\[\sup_{u\in\calT}\int_{0}^{1}\int_{0}^{1}
    (u(x)-u(y))^2K(x-y)\,dxdy<+\infty,\]
with $K(x-y)$ defined in \eqref{eq:K}. Then $\calT$ is relatively compact in $L^2([0,1])$.
\end{lem}

Indeed, one may follow the same proof with $|x-y|^{-1-2s}$ replaced by $K(x-y)$. 

We therefore deduce the convergence of minimizers.

\begin{cor}
Suppose $\calE_\eps(u_\eps)$ is uniformly bounded for a sequence of $\eps\to0^+$. Then there exists a convergent subsequence
\[u_\eps\to{u_0}:=\chi_E-\chi_{E^c}\quad\textin{L^2([0,1])},\]
for some suitable set $E$.

Moreover, if $u_\eps$ minimizes $\calE_\eps$, then $u_0$ minimizes $\calE$.
\end{cor}

\begin{remark}
When $ s\in\left[\frac{1}{2},1\right)$, as in \cite{Savin-Valdinoci12}, up to some multipliers, the energy functional $ {\mathcal H} $ convergences to $\int_0^1 |Du |$ for all $u \in X$. In this case the minimizers coincide with the case of $s=1$ \cite{Ren-Wei00}.
\end{remark}

\section{Computations for the fractional norm $\calH$}
\label{sec:H}

In this section, we compute the derivatives of the $H^s$ norm. The following computations actually hold for kernels more general than the one given in \eqref{eq:K}, as long as
\begin{equation}\label{eq:Kcond}
K(x)=K(-x)=K(1-x)
\end{equation}

\begin{prop}\label{prop:D2H}
Let $N\geq2$ be even and $1\leq{i,j}\leq{N}$ with $i\neq{j}$. Suppose that the kernel $K$ satisfies \eqref{eq:Kcond}. Then
\[\p_{x_ix_j}\calH(U_N)=4(-1)^{i+j-1}K(x_i-x_j)\]
and
\[\p_{x_ix_i}\calH(U_N)=
4\displaystyle\sum_{k=1}^{N-1}(-1)^{k}
    K\left(\dfrac{k}{N}\right).
\]
In particular, the Hessian $D^2\calH(U_N)$ with respect to $(x_1,\dots,x_N)$ is a circulant matrix.
\end{prop}

\begin{proof}
We consider the truncation of the kernel
\[K_M(x-y)=\max\set{K(x-y),M},\]
which is needed for the cancellations of singular terms. Clearly,
\[\lim_{M\to\infty}K_M(x-y)=K(x-y)\]
for a.e. $x,y\in(0,1)$. By Lebesgue Dominated Convergence Theorem (justified by the order $1+2s<2$),
\[\calH(u)
=\lim_{M\to\infty}\calH_M(u)
:=\lim_{M\to\infty}\dfrac14\int_0^1\int_0^1
    (u(x)-u(y))^2K_M(x-y)\,dxdy.\]
Thus it suffices to establish the assertion with the truncated kernel $K_M$.

Since $u^2=1$ a.e., one readily expands
\[\begin{split}
&\quad\;\dfrac14|u(x)-u(y)|^2\\
&=\dfrac{1-u(x)u(y)}{2}\\
&=-\sum_{k=1}^{N}(-1)^{k}H(x-x_k)
    -\sum_{\ell=1}^{N}(-1)^{\ell}H(y-x_\ell)
    -2\sum_{k,\ell=1}^{N}(-1)^{k+\ell}H(x-x_k)H(y-x_\ell).\\
\end{split}\]
As distributions,
\[\begin{split}
&\quad\;
    \p_{x_i}\left(\dfrac12|u(x)-u(y)|^2\right)\\
&=(-1)^{i}\delta(x-x_i)
    +(-1)^{i}\delta(y-x_k)\\
&\quad\;
    +2\sum_{k\neq{i}}(-1)^{k+i}H(x-x_k)\delta(y-x_i)
    +2\sum_{\ell\neq{i}}(-1)^{i+\ell}\delta(x-x_i)H(y-x_\ell)\\
&\quad\;
    +2H(x-x_i)\delta(y-x_i)
    +2\delta(x-x_i)H(y-x_i).
\end{split}\]
Hence, we have
\[\begin{split}
\calH_M(u)
&=-2\sum_{k=1}^{N}(-1)^{k}
        \int_0^1\int_{x_k}^1K_M(x-y)\,dxdy
    -\sum_{k,\ell=1}^{N}(-1)^{k+\ell}
        \int_{x_k}^{1}\int_{x_\ell}^{1}K_M(x-y)\,dxdy,\\
\end{split}\]
as well as the derivatives
\[\begin{split}
\p_{x_i}\calH_M(u)
&=2(-1)^{i}\int_0^1K_M(x-x_i)\,dx
    +4\sum_{k=1}^{N}(-1)^{k+i}\int_{x_k}^{1}K_M(x-x_i)\,dx,\\
\p_{x_ix_j}\calH_M(u)
&=4(-1)^{i+j-1}K_M(x_i-x_j),\\
\p_{x_ix_i}\calH_M(u)
&=2(-1)^{i+1}\int_0^1K_M'(x-x_i)\,dx
    +4\sum_{k=1}^{N}(-1)^{k+i+1}
        \int_{x_k}^{1}K_M'(x-x_i)\,dx
    -4K_M(0)\\
&=2(-1)^{i+1}\left(K_M(1-x_i)-K_M(-x_i)\right)\\
&\quad\;
    +4\sum_{k=1}^{N}(-1)^{k+i+1}\left(K_M(1-x_i)-K_M(x_k-x_i)\right)
    -4K_M(0)\\
&=(-1)^{i}\left(
        2(1-(-1)^{N})K_M(x_i)
        +4\sum_{k\neq{i}}(-1)^{k}K_M(x_k-x_i)
    \right)\\
&=
4\displaystyle\sum_{k\in\set{1,\dots,N}\setminus\set{i}}
    (-1)^{k+i}K_M(x_k-x_i).
\end{split}\]
It suffices to prove that, for the equi-distributed $U_N$, the last expression is independent of $i$. Indeed, 
\[\begin{split}
\p_{x_{2}x_{2}}\calH_M(U_N)
&=4\sum_{k\in\set{1,\dots,N}\setminus\set{2}}
    (-1)^{k+2}K_M(x_k-x_2)\\
&=-4K_M\left(\dfrac{1}{N}\right)
    +4\sum_{k=3}^{N}
        (-1)^{k+2}K_M(x_k-x_2)\\
&=4\sum_{k=1}^{N-2}
    (-1)^{k}K_M(x_{k+2}-x_2)
    -4K_M\left(\dfrac{N-1}{N}\right)\\
&=4\sum_{k=1}^{N-2}
    (-1)^{k}K_M(x_{k+1}-x_1)
    +(-1)^{N-1}4K_M(x_N-x_1)\\
&=4\sum_{k=1}^{N-1}
    (-1)^{k}K_M(x_{k+1}-x_1)\\
&=4\sum_{k=1}^{N-1}
    (-1)^{k}K_M\left(\dfrac{k}{N}\right),
\end{split}\]
by using the properties of $K$ (shared by $K_M$) and the fact that $x_k=\frac{2k-1}{2N}$.
Clearly, one may repeat this argument, shifting more terms to the end, to see that the other diagonal entries $\p_{x_ix_i}\calH_M(U_N)$ is a constant regardless of the value of $i$. The proof is then completed by taking $M\to\infty$.
\end{proof}

\section{Computations related to $\calK$}
\label{sec:K}

\subsection{The Green function}

In this section, we may take any $m\in(-1,1)$. Let $N\geq2$ be even. For $u\in\AN$, we also need to compute the derivatives of $\calK(u)$. Write $v=\Dnh(u-m)$, the solution of
\begin{equation}\label{eq:v}\begin{cases}
-\gamma^2v=u-m
    &\textin(0,1),\\
v(0)=v(1),\\
v'(0)=v'(1),\\
\displaystyle\int_0^1v=0.
\end{cases}\end{equation}
Let $G(x-y)$ be the Green function for the above equation.

\begin{lem}
We have $G(x-y)=\frac{1}{2\gamma^2}B_2 (|x-y|)$, the second Bernoulli polynomial extended periodically hence evenly. More explicitly,
\[G(x-y)=\dfrac{1}{2\gamma^2}\left((x-y)^2-|x-y|+\dfrac16\right).\]
\end{lem}

\begin{proof}
Recall the Fourier series expansion for a function
\[v(x)=\dfrac{v_0}{2}
    +\sum_{k=1}^{\infty}
        (v_k\cos(2\pi{k}x)
        +\tilde{v}_k\sin(2\pi{k}x)),
\]
where, for $k=1,2,\dots$,
\[\begin{split}
v_0&=2\int_0^1v(x)\,dx,\\
v_k&=2\int_0^1v(x)\cos(2\pi{k}x)\,dx,\\
\tilde{v}_k&=2\int_0^1v(x)\sin(2\pi{k}x)\,dx.
\end{split}\]
By comparing the Fourier coefficients of $v$ and $f=u-m$, we have
\[\begin{split}
f_0&=0,\\
f_k&=(2\pi\gamma{k})^2v_k,\\
\tilde{f}_k&=(2\pi\gamma{k})^2\tilde{v}_k.
\end{split}\]
The first equation is satisfied as $m=\int_0^1{u}$. Hence, for $v$ satisfying $\int_0^1v=0$, the unique solution is given by
\[\begin{split}
v(x)
&=\sum_{k=1}^{\infty}
    \left(
        \dfrac{f_k}{(2\pi\gamma{k})^2}\cos(2\pi{k}x)
        +\dfrac{\tilde{f}_k}{(2\pi\gamma{k})^2}\sin(2\pi{k}x)
    \right)\\
&=\sum_{k=1}^{\infty}\dfrac{2}{(2\pi\gamma{k})^2}
    \left(
        \int_0^1f(y)\cos(2\pi{k}y)\,dy\cdot\cos(2\pi{k}x)
        +\int_0^1f(y)\sin(2\pi{k}y)\,dy\cdot\sin(2\pi{k}x)
    \right)\\
&=\int_0^1f(y)
    \left(
        \sum_{k=1}^{\infty}
            \dfrac{2\cos(2\pi{k}(x-y))}{(2\pi\gamma{k})^2}
    \right)
    \,dy.\\
\end{split}\]
Therefore, according to \cite{Abramowitz-Stegun64}, the Green function is
\[\begin{split}
G(x-y)
&=\sum_{k=1}^{\infty}
    \dfrac{2\cos(2\pi{k}(x-y))}{(2\pi\gamma{k})^2}\\
&=\dfrac{1}{4\pi^2\gamma^2}\sum_{k\in\Z\setminus\set{0}}
    \dfrac{e^{2\pi{i}k(x-y)}}{k^2}\\
&=\dfrac{1}{2\gamma^2}B_2(|x-y|)\\
&=\dfrac{1}{2\gamma^2}\left((x-y)^2-|x-y|+\dfrac16\right).
\end{split}\]
\end{proof}

\begin{cor}\label{cor:Kalt}
For $u\in\AN$ where $N\geq1$, one can write
\[\calK(u)=\dfrac12\int_0^1\int_0^1G(x-y)(u(x)-m)(u(y)-m)\,dxdy.\]
Moreover, we have
\[\calK(u)=-\dfrac12\int_0^1\int_0^1(u(x)-u(y))^2G(x-y)\,dxdy
.\]
\end{cor}

\begin{proof}
Indeed,
\[\begin{split}
\calK(u)
&=\dfrac12\int_0^1\abs{\Dnh(u-m)(x)}^2\,dx\\
&=\dfrac12\int_0^1(u(x)-m)(-\gamma^2\Delta)^{-1}(u-m)(x)\,dx\\
&=\dfrac12\int_0^1(u(x)-m)
    \left(\int_0^1G(x-y)(u(y)-m)\,dy\right)\,dx\\
&=\dfrac12\int_0^1\int_0^1G(x-y)(u(x)-m)(u(y)-m)\,dxdy.
\end{split}\]
On the other hand, we have
\[\begin{split}
&\quad\;(u(x)-m)(u(y)-m)\\
&=u(x)u(y)-mu(x)-mu(y)+m^2\\
&=-\dfrac{(u(x)-u(y))^2}{2}+1-m(u(x)-m+m)-m(u(y)-m+m)+m^2\\
&=-\dfrac{(u(x)-u(y))^2}{2}+1-m^2+m(u(x)-m)-m(u(y)-m)
\end{split}\]
so that the second assertion follows from the facts that
\[\int_0^1\int_0^1G(x-y)(u(x)-m)\,dxdy=\int_0^1v(y)\,dy=0\]
and
\[\int_0^1\int_0^1G(x-y)\,dxdy=\int_0^1\Dn(1)\,dy=0.\]
\end{proof}

\subsection{The Hessian}

Following the arguments in \cite{Ren-Wei00} with $G=B_2$, we see that for $u\in\AN$,
\[\p_{x_i}\calK(u)=(-1)^{i-1}2v(x_i).\]

\begin{prop}
The Hessian for $\calK(u)$ with respect to $(x_1,\dots,x_N)$ is given by
\[\p_{x_ix_j}\calK(u)=(-1)^{i-1}2\p_{x_j}v(x_i),\]
with
\[\begin{split}
\p_{x_j}v(x_i)
&=(-1)^{j-1}2G(x_i-x_j)
    +\delta_{ij}\left((1-(-1)^{N})G(x_i)+2\sum_{k=1}^{N}(-1)^{k}G(x_i-x_k)\right).
\end{split}\]
\end{prop}

\begin{proof}
Using the integral representation, we have
\[\begin{split}
\p_{x_j}v(x_i)
&=\p_{x_j}\int_0^1(u(y)-m)G(x_i-y)\,dy\\
&=\p_{x_j}\sum_{k=0}^{N}\int_{x_k}^{x_{k+1}}((-1)^k-m)G(x_k-y)\,dy\\
&=(-1)^{j-1}2G(x_i-x_j)
    +\delta_{ij}\int_0^1(u(y)-m)G'(x_i-y)\,dy\\
&=(-1)^{j-1}2G(x_i-x_j)
    +\delta_{ij}\sum_{k=0}^{N}\int_{x_k}^{x_{k+1}}((-1)^{k}-m)G'(x_i-y)\,dy\\
&=(-1)^{j-1}2G(x_i-x_j)
    +\delta_{ij}\sum_{k=0}^{N}((-1)^k-m)(G(x_i-x_k)-G(x_i-x_{k+1}))\\
&=(-1)^{j-1}2G(x_i-x_j)
    +\delta_{ij}\left((1-(-1)^{N})G(x_i)+2\sum_{k=1}^{N}(-1)^{k}G(x_i-x_k)\right).
\end{split}\]
\end{proof}

\begin{remark}
In view of Corollary \ref{cor:Kalt}, one may also obtain the Hessian of $\calK(u)$ following the proof of Proposition \ref{prop:D2H} with $K=G$ (without the need of a truncation).
\end{remark}

%
%
%
%
%
%

\section{Local minimizers of $\calE$}
\label{sec:localmin}

We aim to prove that for any even $N\geq2$, the equi-distributed $U_N$ is a local minimizer of $\calE$ considered as a functional of $(x_1,\dots,x_N)$, i.e. $D\calE(U_N)=0$, and $D^2\calE(U_N)\geq0$ in the sense of positive semi-definiteness.

\begin{lem}
Let $m=0$. For any even $N\geq2$, $U_N$ is a critical point of $\calE$, i.e.
\[D_{(x_1,\dots,x_N)}\calE(U_N)=0.\]
\end{lem}

\begin{proof}
By the method of Lagrange multiplier, it suffices to verify
\begin{equation}\label{eq:Lagrange}
\p_{x_i}\calH(U_N)
    +\p_{x_i}\calK(U_N)
    +\lambda\p_{x_i}\int_0^1U_N=0,\quad
        i=1,\dots,N,
\end{equation}
for $u=U_N$, $v=\Dn{U_N}$ and some Lagrange multiplier $\lambda\in\R$. Indeed, we have
\[\p_{x_i}\int_0^1U_N
    =\p_{x_i}\int_0^1\sum_{k=0}^{N}(-1)^{k}(H(x-x_k)-H(x-x_{k+1}))\,dx
    =2(-1)^{i-1}\]
and, since
\[\begin{split}
\p_{x_i}\dfrac{(u(x)-u(y))^2}{4}
&=\p_{x_i}\dfrac{1-u(x)u(y)}{2}\\
&=-\dfrac{u(x)\p_{x_i}u(y)+u(y)\p_{x_i}u(x)}{2}\\
&=(-1)^{i}\left(u(x)\delta(y-x_i)+u(y)\delta(x-x_i)\right),
\end{split}\]
we obtain (using the expression in Corollary \ref{cor:Kalt})
\[\begin{split}
&\quad\;\p_{x_i}(\calH(U_N)+\calK(U_N))\\
&=\int_0^1\int_0^1(-1)^{i}
    \left(U_N(x)\delta(y-x_i)+U_N(y)\delta(x-x_i)\right)
    \left(K(x-y)-2G(x-y)\right)\,dxdy\\
&=2(-1)^{i}\int_0^1(K-2G)(x-x_i)U_N(x)\,dx.
\end{split}\]
For $i=1,\dots,N-1$,
\[\begin{split}
&\quad\;\int_0^1(K-2G)(x-x_i)U_N(x)\,dx\\
&=\int_{\frac{i}{N}}^{1+\frac{i}{N}}
    (K-2G)\left(x-x_i-\dfrac{i}{N}\right)U_N(x)\,dx\\
&=\int_{\frac{i}{N}}^{1}
    (K-2G)\left(x-x_i-\dfrac{i}{N}\right)U_N(x)\,dx
+\int_{1}^{1+\frac{i}{N}}
    (K-2G)\left(x-x_i-\dfrac{i}{N}\right)U_N(x)\,dx\\
&=\int_{\frac{i}{N}}^{1}
    (K-2G)(x-x_{i+1})U_N(x)\,dx
+\int_{0}^{\frac{i}{N}}(K-2G)(x-x_{i+1}+1)U_N(x+1)\,dx\\
&=\int_{0}^{1}(K-2G)(x-x_{i+1})U_N(x)\,dx.
\end{split}\]
Therefore, \eqref{eq:Lagrange} is verified with \[\lambda=\int_0^1(K-2G)(x-x_1)U_N(x)\,dx.\]
This completes the proof.
\end{proof}

We now prove Theorem \ref{thm:localmin}.

\begin{proof}[Proof of Theorem \ref{thm:localmin}]
Since $N$ is even, we observe that the Hessian $D^2\calE(U_N)$ is circulant and symmetric. Moreover,
\[
D^2\calE(U_N)
=\begin{bmatrix}
a_0 & a_1 & a_2 & \cdots & a_{N-2} & a_{N-1} \\
a_{N-1} & a_0 & a_1 & a_2 & \cdots & a_{N-2} \\
a_{N-2} & a_{N-1} & a_0 & a_1 & a_2 & \cdots \\
\vdots & \vdots & \vdots & \vdots & \vdots & \vdots \\
a_2 & \cdots & a_{N-2} & a_{N-1} & a_0 & a_1 \\
a_1 & a_2 & \cdots & a_{N-2} & a_{N-1} & a_0
\end{bmatrix},
\]
where
\[a_0=4\sum_{k=1}^{N-1}(-1)^{k}(K-2G)\left(\dfrac{k}{N}\right),\]
\[a_k=a_{N-k}=4(-1)^{k-1}(K-2G)\left(\dfrac{k}{N}\right),\quad
    k=1,\dots,N-1.\]
The eigenvalues are given by
\[\lambda_\ell=\sum_{k=0}^{N-1}a_ke^{i\frac{2\pi{k}\ell}{N}},\quad
    \ell=0,1,\dots,N-1,\]
with corresponding normalized eigenvectors
\[E_\ell=\dfrac{1}{\sqrt{N}}
\begin{pmatrix}
1\\
e^{i\frac{2\pi\ell}{N}}\\
e^{i\frac{2\pi\ell\cdot2}{N}}\\
\vdots\\
e^{i\frac{2\pi\ell\cdot(N-1)}{N}}\\
\end{pmatrix}.
\]
In particular, when $\ell=0$, we have
\[\begin{split}
\lambda_0=\sum_{k=0}^{N-1}a_k=0
\end{split}\]
since $a_0=-\sum_{k=1}^{N-1}a_k$. For $\ell=1,\dots,N-1$, the eigenvalues are given by
\[\begin{split}
\lambda_\ell
&=a_0+\sum_{k=1}^{N-1}a_ke^{i\frac{2\pi{k}\ell}{N}}\\
&=-\sum_{k=1}^{N-1}a_k+\dfrac12\sum_{k=1}^{N-1}
    a_k\left(e^{i\frac{2\pi{k}\ell}{N}}
        +e^{-i\frac{2\pi{k}\ell}{N}}\right)\\
&=-\sum_{k=1}^{N-1}a_k\left(1-\cos\left(\dfrac{2\pi{k}\ell}{N}\right)\right)\\
&=8\sum_{k=1}^{N-1}(-1)^{k}(K-2G)\left(\dfrac{k}{N}\right)
    \sin^2\left(\dfrac{\pi{k}\ell}{N}\right).
\end{split}\]

For any $\ell=1,\dots,N-1$, the contribution from $K$ is bounded in absolute value by
\begin{equation*}
\begin{split}
&\quad\;\abs{8\sum_{k=1}^{N-1}
    (-1)^{k}K\left(\dfrac{k}{N}\right)
    \sin^2\left(\dfrac{\pi{k}\ell}{N}\right)}\\
&\leq{C(s)}\sum_{k=1}^{N-1}\left(
    1+\abs{\dfrac{k}{N}}^{-1-2s}
    +\abs{1-\dfrac{k}{N}}^{-1-2s}\right)
    \sin^2\left(\dfrac{\pi{k}\ell}{N}\right)\\
&\leq{C(s)}N^{1+2s}.
\end{split}\end{equation*}

When $\ell=\frac{N}{2}$, the term with $\sin^2$ is $1$ when $k$ is odd and $0$ otherwise. Hence
\begin{equation*}
\begin{split}
16\sum_{k=1}^{N-1}(-1)^{k-1}G\left(\dfrac{k}{N}\right)\sin^2\left(\dfrac{\pi{k}\ell}{N}\right)
&=\dfrac{8}{\gamma^2}
    \sum_{j=1}^{\frac{N}{2}}\left(
        \dfrac{(2j-1)^2}{N^2}
        -\dfrac{2j-1}{N}
        +\dfrac16
    \right)\\
&=\dfrac{8}{\gamma^2}\left(
        \dfrac{N^2-1}{6N}
        -\dfrac{N}{4}
        +\dfrac{N}{12}
    \right)\\
&=-\dfrac{4}{3\gamma^2N},
\end{split}\end{equation*}
which implies
\[\lambda_{\frac{N}{2}}\leq-\dfrac{4}{3\gamma^2N}+C(s)N^{1+2s}<0\]
provided that
\[\gamma<\dfrac{2}{\sqrt{3C(s)}N^{1+s}}.\]

When $1\leq\ell\leq{N-1}$ and $\ell\neq\frac{N}{2}$, using arithmetico-geometric series (treating trigonometric functions as exponentials) or symbolic computations\footnote{It can be checked, for instance, with the Mathematica code:
\begin{verbatim}
FullSimplify[2Sum[(-1)^(k-1)(k^2/(4n^2) - k/(2n) + 1/6)Sin[k \[Pi] l / (2n)]^2,
{k,1,2n-1}],Assumptions->{n\[Element]Integers,l\[Element]Integers}]//TraditionalForm
\end{verbatim}
},
we see that
\begin{equation*}
16\sum_{k=1}^{N-1}(-1)^{k-1}G\left(\dfrac{k}{N}\right)\sin^2\left(\dfrac{\pi{k}\ell}{N}\right)
=\dfrac{1}{\gamma^2N}\tan^2\left(\dfrac{\pi\ell}{N}\right)
    >0.
\end{equation*}
Hence,
\[\lambda_{\ell}\geq\dfrac{1}{\gamma^2N}\tan^2\left(\dfrac{\pi\ell}{N}\right)-C(s)N^{1+2s}>0,\]
whenever
\begin{equation}\label{eq:gamma0}
\gamma
<\dfrac{1}{\sqrt{C(s)}N^{1+s}}\tan\left(\dfrac{\pi}{N}\right)
=:\gamma_0(N,s).
\end{equation}
This follows from the fact that $|\tan(\pi\ell/N)|$ attains its minimum when $\ell=1$ or $N-1$.

In fact, here $C(s)$ can be taken as $100C_{1,s}$. It stays bounded for $s\in(0,\frac12)$ and tends to zero as $s\to0^+$. 

This completes the proof of Theorem \ref{thm:localmin}.
\end{proof}

%

\section{An explicit solution}
\label{sec:v}

\begin{prop}\label{prop:v}
Suppose $v$ satisfies \eqref{eq:v} with $u=U_N$, $m=0$ and $N$ is even. Then for any $x\in[0,1]$, 
\[v'(x)=-\dfrac{1}{\gamma^2}
    \left(x+2\sum_{k=1}^{N}(-1)^{k}(x-x_k)H(x-x_k)\right)
\]
and
\[v(x)=-\dfrac{1}{\gamma^2}
    \left(-\dfrac{1}{8N^2}+\dfrac{1}{2}x^2
    +\sum_{k=1}^{N}(-1)^{k}(x-x_k)^2H(x-x_k)\right).
\]
In particular, $v'(0)=0$.

\end{prop}

\begin{proof}
Without loss of generality, we may assume that $\gamma=1$. Recall the representation \eqref{eq:uH}, i.e.
\[U_N(x)=1+2\sum_{k=1}^{N}(-1)^{k}H(x-x_k)
    \quad\texton(0,1),\]
with $x_k=\frac{2k-1}{2N}$. Since the Heaviside step function has the integral
\[\begin{split}
\int_{-\infty}^{x}H(t)\,dt
&=\begin{cases}
    x,&\textfor{x\geq0}\\
    0,&\textfor{x\leq0}
  \end{cases}\\
&=xH(x),
\end{split}\]
we have,
\[-v'(x)=-v'(0)+x+2\sum_{k=1}^{N}(-1)^{k}(x-x_k)H(x-x_k)
    \quad\texton(0,1).\]
Similarly, integrating once again from $0$ to $x$,
\begin{equation}\label{eq:vexplicit}
-v(x)=-v(0)-v'(0)x+\dfrac{1}{2}x^2
    +\sum_{k=1}^{N}(-1)^{k}(x-x_k)^2H(x-x_k)
        \quad\texton(0,1).
\end{equation}
Putting $x=1$ and using the condition $v(0)=v(1)$, we find
\[v'(0)=\dfrac12+\sum_{k=1}^{N}(-1)^{k}(1-x_k)^2.\]
As
\[\begin{split}
\sum_{k=1}^{N}(-1)^{k}(1-x_k)^2
=-\dfrac12,
\end{split}\]
we have $v'(0)=0$. In order to find $v(0)$, we integrate \eqref{eq:vexplicit} on $[0,1]$ to obtain
\[v(0)=\dfrac16+\sum_{k=1}^{N}(-1)^{k}
    \dfrac{(1-x_k)^3}{3}
    =\dfrac{1}{8N^2}.\]
This completes the proof.
\end{proof}

\section{An energy growth estimate}
\label{sec:growth}

In this section we prove the following energy estimate for the local minimizer $U_N$. Theorem \ref{thm:EUN} is a direct consequence of the following two lemmata. Below we write
$f\sim{g}$ if there exists a constant $C>0$ such that $0<C^{-1}f\leq{g}\leq{C}f$.

\begin{lem}
We have
\[\calK(U_N)=\dfrac{1}{24\gamma^2N^2}.\]
\end{lem}

\begin{proof}
By Proposition \ref{prop:v}, $v$ satisfies the Neumann boundary condition. Hence, the calculations in \cite{Ren-Wei00} give the desired formula.
\end{proof}

\begin{lem}
There exists a constant $C>0$ such that
\[\calH(U_N)\sim{N}^{2s}.\]
\end{lem}

\begin{proof}
For the integral
\[\calH(U_N)=\int_0^1\int_0^1|u(x)-u(y)|^2K(x-y)\,dxdy,\]
we first observe that the non-zero contributions come from the region
\[\calR=\set{(x,y)\in[0,1]^2:u(x)u(y)=-1}.\]
It suffices to consider only the case $u(x)=1$ and $u(y)=-1$, each number of intervals is $\frac{N}{2}$. Thus $\calR$ can be decomposed to a union of $\left(\frac{N}{2}\right)^2$ rectangles, the product of intervals $(x_j,x_{j+1})\times(x_k,x_{k+1})$ where $j$ and $k$ have opposite parity.

By the periodicity of $U_N$ and noting that the far away interactions are comparable to their adjacent ones, all these interactions in $\calR$ can be grouped together and computed by

%
\[\begin{split}
\calH(U_N)
&\sim{N}\int_0^{x_1}\int_{x_1}^{x_{\frac{N}{2}}}
    \left(\dfrac{1}{(x-y)^{1+2s}}+O_s(1)\right)\,dxdy\\
&\sim{N}\int_0^{x_1}\left(\dfrac{1}{(x_3-y)^{2s}}+O_s(1)\right)\,dy\\
&\sim{N}^{2s}.
\end{split}\]
\end{proof}

\appendix
\section{A brief derivation of the free energy}

As a physical motivation, here we include a brief derivation for the free energy with nonlocal diffusion using the density functional theory of Ohta--Kawasaki \cite{Ohta-Kawasaki86}, following closely Choksi--Ren \cite{Choksi-Ren03}. We will use the notations as in \cite{Choksi-Ren03} and point out the notable differences.

Suppose a diblock copolymer consists of chains of monomers $A$ and $B$, and the melt lives on $\Omega\subset\R^3$. Write $N$ to be the index of polymerization. The intervals occupied by the $A$- and $B$-monomers are denoted $\calI_A=[0,N_A]$ and $\calI_B=[N_A,N]$. Write $N_B=N-N_A$. We assume the Kuhn statistical lengths for the A and B monomers are the same and equal $1$. For each $i=1,\dots,n$, a copolymer chain $r_i:[0,N]\to\R^3$ is a continuous function for each $i$. The phase space is then
\[\Gamma=\set{r=(r_1,\dots,r_n):r_i\in C([0,N],\R^3)},\]
equipped with the measure
\[d\mu=\underbrace{(dx\times d\calP_0)\times\cdots\times(dx\times d\calP_0)}_{n},\]
where $d\calP_0$ is, instead of the classical Wiener measure of the Brownian motion, one driven by the isotropic $2s$-stable L\'{e}vy process starting at the origin. As in \cite{Choksi-Ren03},
the Hamiltonian can be written as
\[H(r)=\sum_{k,m\in\set{A,B}}\int_{\Omega}\dfrac{V^{km}}{2\rho_0}\rho_k(x,r)\rho_m(x,r)\,dx\]
where $V^{km}>0$ denote the monomer interaction parameters, $\rho_0=nN/|\Omega|$ denote the average monomer number density, 
and $\rho_k(x,r)=\sum_{i=1}^{n}\int_{\calI_k}\delta(x-r_i(\tau))\,d\tau$ are the microscopic density fields.
In this way, the Gibbs canonical distribution is
\[D(r)=\dfrac{1}{Z}e^{-\beta H(r)}\]
and the free energy of the system is
\[-\beta^{-1}\log{Z},\]
where
\[Z=\int_{\Gamma}e^{-\beta H(r)}\,d\mu,\]
and $\beta$ is the reciprocal of the absolute temperature (measured in the energy unit so the Boltzmann constant equals $1$).

In the self-consistent mean field theory we choose the class of distributions generated by an external fields pair $U=(U^A,U^B)$, acting on monomers $A$ and $B$ respectively and there is no interaction between them. By adding a suitable constant we assume that
\[\sum_{k\in\set{A,B}}\dfrac{N_k}{N}\int_{\Omega}U^k(x)\,dx=0.\]
It induces
\begin{itemize}
\item the Hamiltonian on $\Gamma$,
	\[H_U(r)=\sum_{i=1}^{n}\sum_{k\in{A,B}}\int_{\calI_k}U^k(r_i(\tau))\,d\tau;\]
\item a Gibbs canonical distribution and the corresponding partition function,
	\[D_U(r)=\dfrac{1}{Z_U}e^{-\beta H_U(r)},
		\quad Z_U=\int_{\Gamma}e^{-\beta H_U(r)}\,d\mu;\]
\item the expectations of the microscopic density fields,
	\[\angles{\rho_k(x)}_U=\int_{\Gamma}\rho_k(x,r)D_U(r)\,d\mu;\]
\item the average internal energy under $D_U$,
	\[\angles{H}_U:=\int_{\Gamma}H(r)D_U(r)\,d\mu=\int_{\Omega}\dfrac{V^{km}}{2\rho_0}\angles{\rho_k(x)}_U\angles{\rho_m(x)}_U\,dx.\]
\item the approximate free energy as a functional of $U$ (via the variational principle of $D$),
	\[\begin{split}
	F(U)&=\angles{H}_U-\beta^{-1}S(D_U)\\
	&=\int_{\Omega}\left(\sum_{k,m\in\set{A,B}}\dfrac{V^{km}}{2\rho_0}\angles{\rho_k(x)}_U\angles{\rho_m(x)}_U-\sum_{k\in{A,B}}U^k(x)\angles{\rho_k(x)}_U\right)\,dx
		-\beta^{-1}\log Z_U,
	\end{split}\]
	where $S(D_U)=-\int D_U\log D_U\,d\mu$ is the entropy of $D_U$.
\end{itemize}

Our goal is to express $F(U)=\angles{H}_U-\beta^{-1}S(D_U)$ in terms of $\angles{\rho_k(x)}_U$. As in the classical model, one has
\[\dfrac{\delta(-S(D_U))}{\delta(\angles{\rho}_U)}=-\beta{U},\]
meaning that for the second entropy term it suffices to express $\beta{U}$ in terms of $\angles{\rho}_U$ and integrate with respect to $\angles{\rho}_U$.

The computation of $F(U)$ is done by the Feynman--Kac integration theory. In terms of the solutions $q_U$ and $q_U^*$ of the backward and forward parabolic equations, one obtains formulae for $Z_U$ and $\angles{\rho_k(x)}_U$. More precisely, let $Q_U(y,\tau,z,t)$ be the fundamental solution of the backward equation
\[(Q_U)_\tau-(-\Delta_y)^{s}Q_U-\beta UQ_U=0,
	\quad Q_U(y,t,z,t)=\delta(y-z),\]
where $U(y,\tau)=U^k(y)$ if $\tau\in\calI_k$, $k\in\set{A,B}$. Then $q_U(y,\tau)=\int_{\Omega}Q_U(y,\tau,z,N)\,dz$ solves the backward equation
\[(q_U)_\tau-(-\Delta)^{s}q_U-\beta{U}q_U=0,
	\quad q_U(y,N)=1,
	\quad (y,\tau)\in\Omega\times(0,N)\]
and $q_U^*(y,\tau)=\int_{\Omega}Q_U(z,0,y,\tau)\,dz$ solves the forward equation
\[(q_U^*)_\tau+(-\Delta)^{s}q_U^*+\beta{U}q_U^*=0
	\quad q_U^*(y,0)=1,
	\quad (y,\tau)\in\Omega\times(0,N).\]
Here we impose the Dirichlet boundary conditions that $q_U$ and $q_U^*$ vanish outside $\Omega$. A probabilistic calculation reveals that
\[Z_U=\left(\int_{\Omega}q_U(y,0)\,dy\right)^{n}
	=\left(\int_{\Omega}q_U^*(y,N)\,dy\right)^{n}\]
and
\[\angles{\rho_k(x)}_U=\dfrac{n}{Z_U^{\frac1n}}\int_{\calI_k}q_U(x,\tau)q_U^*(x,\tau)\,d\tau.\]

Let us apply the first approximation, namely linearizing the dependence of $U$ around $0$,
\[\begin{split}
&\quad\,\angles{\rho_k(x)}_U\\
&\approx
	\angles{\rho_k(x)}_0
	+\left.\dfrac{d}{d\eps}\right|_{\eps=0}\angles{\rho_k(x)}_{0+\eps{U}}\\
&=\angles{\rho_k(x)}_0
	+\dfrac{\angles{\rho_k(x)}_0}{n}\int_{\Omega}\sum_{m\in\set{A,B}}\angles{\rho_m(y)}_0\beta{U}^m(y)\,dy
	+\dfrac{n}{Z_0^{\frac1n}}\int_{\calI_k}(pq_0^*+p^*q_0)(x,\tau)\,d\tau,
\end{split}\]
using the same calculations (except changing the fundamental solution to the fractional one) of \cite{Choksi-Ren03}, where $p$ and $p^*$ solve respectively
\[p_\tau-(-\Delta)^{s}p-\beta{U}q_0=0,
	\quad p(x,N)=0,\]
\[p_\tau^*+(-\Delta)^{s}p^*+\beta{U}q_0=0,
	\quad p^*(x,0)=0,\]
and can be expressed as
\[p(x,\tau)=-\beta\int_\tau^{N}\int_{\Omega}Q_0(x,\tau,y,t)q_0(y,t)U(y,t)\,dydt,\]
\[p^*(x,\tau)=-\beta\int_0^{\tau}\int_{\Omega}Q_0(y,t,x,\tau)q_0^*(y,t)U(y,t)\,dydt.\]

In the second approximation one takes the thermodynamic limit, letting $\Omega\to\R^3$, $n\to\infty$, while keeping $\frac{n}{|\Omega|}=\frac{\rho_0}{N}$ unchanged. Then
\[Q_0(y,\tau,z,N)\to\calK_s(y-z,\tau-t),\]
the fractional heat kernel in $\R^3$, whose Fourier transform in $\R^3$ equals
\[\hat\calK_s(\xi)=e^{-t(2\pi|\xi|)^{2s}}.\]
We also have
\[q_0\to1,\quad 
	q_0^*\to 1,\quad
	\dfrac{n}{Z_0^{\frac1n}}\to\dfrac{\rho}{N},\quad
	\angles{\rho_k(x)}_0\to\dfrac{N_k}{N}\rho_0=:\bar\rho_k,\]
\[\sum_{m\in\set{A,B}}\dfrac{\angles{\rho_k(x)}_0}{n}\int_{\Omega}\angles{\rho_m(y)}_0\beta U^m(y)\,dy\to0,\]
\[p(x,\tau)\to -\int_\tau^N\left(\calK_s(\cdot,\tau-t)\ast \beta U(\cdot,t)\right)(x)\,dt,\]
\[p^*(x,\tau)\to -\int_0^\tau\left(\calK_s(\cdot,\tau-t)\ast \beta U(\cdot,t)\right)(x)\,dt,\]
and, more importantly,
\[\angles{\rho_k}_{U}\approx \bar\rho_k-\dfrac{\rho_0}{N}\sum_{m\in{A,B}}R_{km}\ast(\beta U^m),\]
where
\[R_{km}(z):=\int_{\calI_k}\int_{\calI_m}\calK_s(z,\tau-t)\,dtd\tau.\]
whose Fourier transform is
\[\hat{R}_{km}(\xi)=
\begin{cases}
2(2\pi|\xi|)^{-4s}h\left((2\pi|\xi|)^{2s}N_k\right),
	&\textif k=m,\\
(2\pi|\xi|)^{-4s}g\left((2\pi|\xi|)^{2s}N_k,(2\pi|\xi|)^{2s}N_m\right)
	&\textif k\neq{m}.
\end{cases}
\]
Here $h(s_1)=e^{-s_1}+s_1-1$ and $g(s_1,s_2)=(1-e^{-s_1})(1-e^{-s_2})$.

From this point, one may apply the third approximation, namely the long and short wave expansions, to obtain that
\[h(s_1)\approx s_1,
	\quad
	g(s_1,s_1)\approx 1
	\quad 
	\textif s_1,s_2\gg 1,\]
\[h(s_1)\approx \dfrac{s_1^2}{2}-\dfrac{s_1^3}{6},
	\quad
	g(s_1,s_2)\approx (s_1-\dfrac{s_1^2}{2})(s_2-\dfrac{s_2^2}{2})
	\quad
	\textif s_1,s_2\ll 1.\]
Thus
\[
\hat{T}(\xi):=(\hat{R})^{-1}(\xi)
	\approx
	\dfrac{(2\pi|\xi|)^{2s}}{N}K+\dfrac{1}{(2\pi|\xi|)^{2s}N^3}L,\]
where
\[K=\dfrac12\begin{pmatrix}a^{-1} & 0 \\ 0 & b^{-1}\end{pmatrix},
\quad
L=\dfrac32\begin{pmatrix}a^{-2} & -(ab)^{-1} \\ -(ab)^{-1} & b^{-2}\end{pmatrix},
\quad
a=\dfrac{N_A}{N},
\quad
b=\dfrac{N_B}{N}.
\]

By taking the inverse Fourier transform,
\[T=\dfrac{1}{N}(-\Delta)^{s}K+\dfrac{1}{N^3}(-\Delta)^{-s}L.\]
Since
\[\beta U^k(x)
	\approx
		-\dfrac{N}{\rho_0}\sum_{m\in\set{A,B}}T^{km}\left(\angles{\rho_m}_U-\bar\rho_m\right)(x),\]
integrating it yields
\begin{multline*}
-S(D_U)+S(D_0)
\approx
	\dfrac{1}{2\rho_0}\int_{\R^3}
	\Bigg(
		\sum_{k\in\set{A,B}}K^{kk}\angles{\rho_k}_U(-\Delta)^{s}\angles{\rho_k}_U
\\	
		+\sum_{k,m\in\set{A,B}}\dfrac{L^{km}}{N^2}
		\left(
			\angles{\rho_k}_U-\bar\rho_k
		\right)
			(-\Delta)^{-s}
			\left(
				\angles{\rho_m}_U-\bar\rho_m
			\right)
	\Bigg).
\end{multline*}
The constant $S(D_0)$ can be dropped. This finally gives rise to the Ohta--Kawasaki free energy with the fractional gradient, upon using the incompressibility constraint.

We remark that the $H^{-s}$-norm behaves in a similar way as the $H^{-1}$-norm which has been used classically, which we have decided to use for the simplicity of the mathematical treatment.

\section*{Acknowledgement}
\thanks{All authors are partially supported by NSERC-2018-03773. HC thanks his advisors Prof. Juncheng Wei and Prof. Nassif Ghoussoub for their constant encouragement and support. We thank the anonymous referees for carefully reading the manuscript and giving valuable comments.}

 \newcommand{\noop}[1]{}

\end{document}